\documentclass[reqno]{amsart}
\usepackage{float}
\usepackage{amsmath}
\usepackage{graphicx}
\usepackage{latexsym}
\usepackage{amsfonts}
\usepackage{amssymb}
\theoremstyle{plain}
\newtheorem{theorem}{Theorem}

\newtheorem{lemma}[theorem]{Lemma}
\newtheorem{proposition}[theorem]{Proposition}
\theoremstyle{definition}

\newtheorem{remark}[theorem]{Remark}
\newtheorem*{remark*}{Remark}

\begin{document}
\title[Transition phenomena for ladder epochs]{Transition phenomena for ladder epochs
of random walks with small negative drift}
\author[Wachtel]{Vitali Wachtel}
\address{Mathematical Institute, University of Munich,
Theresienstrasse 39, D-80333, Munich, Germany}
\email{wachtel@math.lmu.de}
\date{\today }

\begin{abstract}
For a family of random walks $\{S^{(a)}\}$ satisfying $\mathbf{E}S_1^{(a)}=-a<0$ 
we consider ladder epochs $\tau^{(a)}=\min\{k\geq1: S_k^{(a)}<0\}$. We study the 
asymptotic, as $a\to0$, behaviour of $\mathbf{P}(\tau^{(a)}>n)$ in the case when $n=n(a)\to\infty$.
As a consequence we obtain also the growth rates of the moments of $\tau^{(a)}$.
\end{abstract}

\maketitle

\section{Introduction and statement of results}
\subsection{Background and purpose}
Let $X,X_1,X_2,\ldots$ be independent identically distributed random variables. 
Let $S=\{S_n,\,n\geq0\}$ denote the random walk with increments $X_i$, that is,
$$
S_0:=0,\ S_n:=\sum_{i=1}^nX_i.
$$
Let us first recall what is known on the first descending ladder epoch $\tau$ of $S$, i.e.,
\begin{equation}
\label{tau}
\tau:=\min\{k\geq1: S_k<0\}.
\end{equation}
It is well-known (see, for example, \cite[Theorem 17.1]{Sp01}) that
\begin{equation*}
\mathbf{P}(\tau<\infty)=1 \Leftrightarrow\sum_{k=1}^\infty k^{-1}\mathbf{P}(S_k<0)=\infty.
\end{equation*}
Under the latter condition Rogozin \cite{Rog71} has studied the asymptotic, as $n\to\infty$, behaviour of the tail probability
$\mathbf{P}(\tau>n)$. In particular, 
\begin{equation}
\label{Rogozin}
\lim_{n\to\infty}\frac{1}{n}\sum_{k=1}^n\mathbf{P}(S_k\geq0)=\rho\in(0,1]
\Leftrightarrow
\mathbf{P}(\tau>n)=n^{\rho-1}\ell(n),
\end{equation}
where $\ell$ is slowly varying at infinity.
Also,
$\lim_{n\to\infty}\frac{1}{n}\sum_{k=1}^n\mathbf{P}(S_k\geq0)=0$ is equivalent to
the relative stability of $\tau$. The latter means that the function $x\mapsto\int_0^x\mathbf{P}(\tau>u)du$
is slowly varying at infinity. But this statement does not give any information on the asymptotic behaviour of $\mathbf{P}(\tau>n)$
in this case.

The situation when $\mathbf{E}\tau<\infty$, which is a particular case of the relative stability,
was considered by Embrechts and Hawkes \cite{EH82}. There it has been shown that 
$$
\mathbf{P}(\tau>n)\sim n^{-1}\mathbf{P}(S_n>0)\exp\Bigl\{\sum_{j=1}^\infty r^j j^{-1}\mathbf{P}(S_j\geq0)\Bigr\},
$$
under certain conditions on the sequence $\{\mathbf{P}(S_n>0),n\geq1\}$.

If the expectation $\mathbf{E}X$ is finite, then the condition $\sum_{k=1}^\infty k^{-1}\mathbf{P}(S_k\leq0)=\infty$
is equivalent to the inequality $\mathbf{E}X\leq0$, see again \cite[Theorem 17.1]{Sp01}. 
If $\mathbf{E}X=0$ and $X$ belongs to the domain of attraction of a stable law, then $\lim_{n\to\infty}\mathbf{P}(S_n\geq0)\in(0,1)$. 
This yields that $\lim_{n\to\infty}\frac{1}{n}\sum_{k=1}^n\mathbf{P}(S_k\geq0)=\rho\in(0,1)$. 
Then, using (\ref{Rogozin}), we conclude that 
\begin{equation}
\label{zeromean}
\mathbf{P}(\tau>n)=n^{\rho-1}\ell(n).
\end{equation}
If $\mathbf{E}X<0$, then $\mathbf{E}\tau$ is finite, see \cite[Proposition 18.1]{Sp01}. 
In this case of negative drift, Doney \cite{Don89} has applied the results from \cite{EH82} to two special classes of random walks:
He has shown that if $\mathbf{E}X\in(-\infty,0)$ and $\mathbf{P}(X>x)$ 
is regularly varying at infinity with index $\alpha<-1$, then, as $n\to\infty$,
\begin{equation}
\label{Don1}
\mathbf{P}(\tau>n)\sim \mathbf{E}\tau\mathbf{P}(X>-n\mathbf{E} X)\quad\text{as }n\to\infty.
\end{equation}
Besides this case of regularly varying tail, Doney found the asymptotics of $\mathbf{P}(\tau>n)$
for random walks having negative drift and satisfying the following condition: If the equation
$\frac{d}{dh}\mathbf{E}e^{hX}=0$ has a positive solution, say $h_0$, then
\begin{equation}
\label{Don2}
\mathbf{P}(\tau>n)\sim C\Bigl(\frac{\mathbf{E}\mu^\tau-1}{\mu-1}\Bigr)
\mu^{-n}n^{-3/2}\quad\text{as }n\to\infty,
\end{equation}
where $\mu=1/\mathbf{E}e^{h_0X}$ and $C$ is a constant depending on $\mathbf{E}e^{hX}$.
The latter relation was generalised by Bertoin and Doney \cite{BD96} to the case when
$\frac{d}{dh}\mathbf{E}e^{hX}<0$ for all $h>0$ such that $\mathbf{E}e^{hX}<\infty$.

It should be noted that \cite{Don89} and \cite{BD96} are devoted to the study of the asymptotic behaviour of $\mathbf{P}(\tau_x>n)$
for any fixed $x\geq0$, where $\tau_x:=\min\{k\geq1: S_k<-x\}.$ The main result can be stated
as follows: If $X$ satisfies the conditions stated before (\ref{Don1}) or (\ref{Don2}), 
then there exists a function $U$ such that
\begin{equation*}
\lim_{n\to\infty}\frac{\mathbf{P}(\tau_x>n)}{\mathbf{P}(\tau>n)}\to U(x).
\end{equation*}

Studying the asymptotic, as $n\to\infty$, behaviour of $\mathbf{P}(\tau>n)$, one hopes to get a 
good approximation for large but finite values of $n$. The quality of such approximation depends on different
parameters of the random walk. It follows from the papers mentioned above that the asymptotic behaviour of $\mathbf{P}(\tau>n)$
depends crucially on whether $\mathbf{E}X=0$ or $<0$. Therefore, it would be very useful to clarify the influence of
$\mathbf{E}X$ on $\mathbf{P}(\tau>n)$ in the case when that expectation is quite small. We illustrate the problem with 
the following concrete example. Let $S$ be a random walk with $\mathbf{E}X=10^{-3}$ and we want to calculate the 
quantity $\mathbf{P}(\tau>10^5)$. Here one has two possibilities: On the one hand, one can say that the 
expectation is so small, that we may apply asymptotic relations for zero mean random walks. And on the other hand, 
we can say that the expectation is negative and we should use formulas (\ref{Don1}) or (\ref{Don2}), depending on the 
tail behaviour of $X$. But how to decide, which approximation is better for these values of $\mathbf{E}X$ and $n$? 
This question leads to the following mathematical problem: What can be said on the asymptotic behaviour of 
$\mathbf{P}(\tau>n)$ in the case when $\mathbf{E}X\to0$ and $n\to\infty$ simultaneously? 

In the present paper we consider this problem in the case when the random walk's increment belongs to the domain of attraction of a stable law.
We shall show that there exists a function $f$ such that
\begin{itemize}
\item[(a)] if $n\ll f(\mathbf{E}X)$, then one has to use (\ref{zeromean}),
\item[(b)] if $n\gg f(\mathbf{E}X)$, then one has to use formulas for random walks with negative drift,
\item[(c)] if $n\sim vf(\mathbf{E}X)$, $v\in(0,\infty)$, then one has to use (\ref{zeromean}), but with a correction
factor depending on $v$.
\end{itemize}
The last point seems to be the most interesting one: It describes {\it transition phenomena} for the ladder epoch $\tau$,
which appear in the case of small drift. 

Our {\it main result}, Theorem~\ref{T.normal}, is devoted to the study of this transition: 
There it will be clarified how the function $f$ and the correction factor look like. As a consequence we will get the claim in (a).
Furthermore, Theorem~\ref{T.normal} allows one to determine the asymptotic, as $\mathbf{E}X\to0$, behaviour of 
some moments $\mathbf{E}\tau^r$, 
see Theorem~\ref{T.exp}. The expectation $\mathbf{E}\tau$ is of particular interest, since it appears in asymptotic relations
connected to the claim in (b), see Theorems~\ref{LD_stable}, \ref{LD_var} and \ref{tail} below.

\subsection{Transition phenomena}
We start with a more precise description of our model of random walks with asymptotically small drift.
We shall consider a family of random walks 
$\mathcal\{S^{(a)},a\in[0,a_0]\}$ with drift $-a$, that is, 
$\mathbf{E}S_1^{(a)}=-a$, and investigate the asymptotic, as $a\to0$, behaviour of the
probability $\mathbf{P}(\tau^{(a)}>n)$ for $n=n(a)$, where $\tau^{(a)}$ is the first descending 
ladder epoch of $S^{(a)}$, as in (\ref{tau}).

Let $X^{(a)}$ denote a random variable, which is distributed as the increments of the random walk $S^{(a)}$.
It is easy to see that if $X^{(a)}$ converges in distribution, as $a\to0$, to $X^{(0)}$, then, for every fixed $n$,
\begin{equation}
\label{fixed}
\mathbf{P}(\tau^{(a)}>n)\sim\mathbf{P}(\tau^{(0)}>n)
\text{ as }a\to0.
\end{equation} 
A more interesting problem consists in investigating the asymptotic behaviour of the tail 
probability $\mathbf{P}(\tau^{(a)}>n)$ when $n=n(a)\to\infty$ as $a\to0$. The answer to this 
question depends on the structure of the family $\mathcal\{S^{(a)},a\in[0,a_0]\}$. 

In this paper we shall assume that 
there exists a random variable $X$ with zero mean such that the random variables
$X^{(a)}$ and $X-a$ have the same distribution for all $a\in[0,a_0]$.
This yields that the random variables $S_n^{(a)}$ and $S_n^{(0)}-na$ are equal in distribution for
all $a\in[0,a_0]$ and $n\geq1$.
Furthermore, we restrict ourselves from now on to so-called {\it asymptotically stable random walks}.
Namely, we shall always assume that the distribution of $X$ belongs to the domain 
of attraction of a stable law with characteristic function
\begin{equation}
\label{std}
G_{\alpha ,\beta }\mathbb{(}t\mathbb{)}:=\exp \left\{ -|t|^{\,\alpha
}\left( 1-i\beta \frac{t}{|t|}\tan \frac{\pi \alpha }{2}\right) \right\}
\end{equation}
with $\alpha\in(1,2]$ and $|\beta|\leq1$.
In this case we write $X\in \mathcal{D}\left( \alpha ,\beta \right)$.

Let $\left\{ c_{n},n\geq 1\right\} $ denote the
sequence of positive integers specified by the relation
\begin{equation}
c_{n}:=\inf \left\{ u\geq 0:u^{-2}V(u)\leq n^{-1}\right\} ,  \label{Defa}
\end{equation}%
where
\begin{equation*}
V(u):=\int_{-u}^{u}x^{2}\mathbf{P}(X\in dx),\ u>0.
\end{equation*}%
It is known (see, for instance, \cite[Ch. XVII, \S 5]{FE}) that the function 
$V$ is regularly varying at infinity with index $2-\alpha$ for every 
$X\in \mathcal{D}(\alpha ,\beta)$. This implies that $\left\{c_{n},n\geq
1\right\} $ is regularly varying with index $\alpha ^{-1}$, i.e.
there exists a function $l_{1},$ slowly varying at infinity, such that
\begin{equation}
c_{n}=n^{1/\alpha }l_{1}(n).  \label{asyma}
\end{equation}
In addition, the scaled sequence $\left\{ S^{(0)}_{n}/c_{n},\,n\geq 1\right\} $
converges in distribution, as $n~\rightarrow~\infty ,$ to the stable law
corresponding to $G_{\alpha ,\beta }$ in (\ref{std}). 

Let $\{Y_{\alpha,\beta}(t),\,t\geq0\}$ denote a stable Levy
process such that $Y_{\alpha,\beta}(1)$ distributed according to (\ref{std}).

It is known, see \cite[Proposition 17.5]{Sp01}, that the generating function of 
the sequence $\{\mathbf{P}(\tau^{(a)}>n),n \geq0\}$ satisfies the identity
\begin{equation}
\label{GenF}
\sum_{n=0}^\infty\mathbf{P}(\tau^{(a)}>n)z^n=
\exp\Bigl\{\sum_{n=1}^\infty\frac{z^n}{n}\mathbf{P}(S^{(a)}_{n}\geq0)\Bigr\},
\quad z\in(0,1).
\end{equation}
Thus, for every $n\geq1$, the probability $\mathbf{P}(\tau^{(a)}>n)$ 
is determined by $\{\mathbf{P}(S^{(a)}_{k}\geq0),1\leq k\leq n\}$.
From the definition of the family $S^{(a)}$ and from the asymptotic stability of $\{ S^{(0)}_{n},n\geq0\}$ we 
conclude that
\begin{equation}
\label{Rho}
\mathbf{P}(S_n^{(a)}\geq0)\sim\mathbf{P}(S_n^{(0)}\geq0)
\sim\mathbf{P}(Y_{\alpha,\beta}(1)\geq0)=:\rho\in(0,1)
\end{equation}
for $n=n(a)\to\infty$ satisfying $na/c_n\to0$. Hence, one can expect that
\begin{equation}
\label{tau0}
\mathbf{P}(\tau^{(a)}>n)\sim\mathbf{P}(\tau^{(0)}>n)=n^{\rho-1}\ell(n),
\end{equation}
where in the second step we have used (\ref{Rogozin}). Furthermore, if
$na/c_n\to u\in(0,\infty)$, then 
$$
\mathbf{P}(S_n^{(a)}\geq0)\sim\mathbf{P}(Y_{\alpha,\beta}(1)\geq u)>0.
$$ 
In this case
one expects, although this conjecture is not as obvious as (\ref{tau0}),
that
\begin{equation}
\label{tau1}
\mathbf{P}(\tau^{(a)}>n)\sim\mathbf{P}(\tau^{(0)}>n)G(u)
\end{equation}
for some function $G$.

The following theorem confirms the conjectures (\ref{tau0}) and (\ref{tau1}). 
\begin{theorem}\label{T.normal}
Suppose $X\in\mathcal{D}(\alpha,\beta)$.
If $n=n(a)$ is such that 
\begin{equation}
\label{condition}
\lim_{a\to0}\frac{an}{c_n}=u\in[0,\infty),
\end{equation}
then
\begin{equation}
\label{T.n1}
\lim_{a\to0}\frac{\mathbf{P}(\tau^{(a)}>n)}{\mathbf{P}(\tau^{(0)}>n)}=(1-F_{\alpha,\beta}(u)),
\end{equation}
where the distribution function $F_{\alpha,\beta}$ can be described by the equality
\begin{align}
\label{T.n22}
\nonumber
\int_0^\infty e^{-\lambda x}x^{\rho-1}(1-F_{\alpha,\beta}(x^{1-1/\alpha}))dx=\hspace{4cm}\\
C\exp\Bigl\{-\int_0^\infty\frac{1-e^{-\lambda t}}{t}\mathbf{P}(Y_{\alpha,\beta}(t)-t>0)dt\Bigr\},\ \lambda\geq0
\end{align}
with $\rho$ defined in (\ref{Rho}) and with $C$ specified by the condition $F_{\alpha,\beta}(0)=0$.
\end{theorem}
The existence of the limit in (\ref{T.n1}) is an easy consequence of the invariance principle 
for random walks conditioned to stay positive, which was proved by Doney \cite{Don85}. The most difficult
part of the proof is the derivation of characterisation (\ref{T.n22}) of the limiting distribution $F_{\alpha,\beta}$,
see Section~\ref{s_T.normal}.

It follows from (\ref{asyma}) that (\ref{condition}) is equivalent to
$$
n\sim u^{\alpha/(\alpha-1)}\Bigl(\frac{1}{a}\Bigr)^{\alpha/(\alpha-1)}l^*\Bigl(\frac{1}{a}\Bigr)
\quad\text{as } a\to0,
$$
where $l^*$ is slowly varying at infinity, which is determined by $l_1$. Therefore, the statement of Theorem~\ref{T.normal} can 
be reformulated as follows: If $n=n(a)$ satisfies
\begin{equation}
\label{con}
n\sim v\Bigl(\frac{1}{a}\Bigr)^{\alpha/(\alpha-1)}l^*\Bigl(\frac{1}{a}\Bigr)
\quad\text{as }a\to0
\end{equation}
for some $v\geq0$, then
\begin{equation}
\label{corr}
\lim_{a\to0}\frac{\mathbf{P}(\tau^{(a)}>n)}{\mathbf{P}(\tau^{(0)}>n)}=
\Bigl(1-F_{\alpha,\beta}\bigl(v^{1-1/\alpha}\bigr)\Bigr).
\end{equation}
In particular, if (\ref{con}) holds with $v=0$, then $\mathbf{P}(\tau^{(a)}>n)\sim\mathbf{P}(\tau^{(0)}>n)$.
Roughly speaking, (\ref{zeromean}) give a rather good approximation in the case when $n$ is much smaller than
$\Bigl(\frac{1}{a}\Bigr)^{\alpha/(\alpha-1)}l^*\Bigl(\frac{1}{a}\Bigr)$. But if 
$\Bigl(\frac{1}{a}\Bigr)^{\alpha/(\alpha-1)}l^*\Bigl(\frac{1}{a}\Bigr)$ and $n$ are comparable, then one has to use a correction factor, given
by the right hand side of (\ref{corr}). To calculate this correction for concrete values of $v$ one has to know the form of the 
distribution function $F_{\alpha,\beta}$.
We are able to give an explicit expression for $F_{\alpha,\beta}$ only in some special cases: We shall see in the proof of Theorem~\ref{T.normal} that 
$$
1-F_{\alpha,\beta}(u)=\mathbf{P}\Bigl(\inf_{t\leq1}(M_{\alpha,\beta}(t)-ut)\geq0\Bigr),
$$
where $\{M_{\alpha,\beta}(t),t\in[0,1]\}$ is the meander of $Y_{\alpha,\beta}$.
Using the construction of the meander via the limit of conditioned distributions of 
the original process $Y_{\alpha,\beta}$,
we shall show that
$$
1-F_{2,0}(u)= u\int_u^\infty v^{-2}e^{-v^2/2}dv
$$
and
$$
1-F_{\alpha,1}(u)=\frac{u^{1/(\alpha-1)}}{(\alpha-1)g_{\alpha,1}(0)}\int_u^\infty v^{-\alpha/(\alpha-1)}g_{\alpha,1}(v)dv,
\quad \alpha\in(1,2),
$$
where $g_{\alpha,\beta}$ denotes the density function of the random variable $Y_{\alpha,\beta}(1)$.
For all other values of $\alpha$ and $\beta$ the explicit form of $F_{\alpha,\beta}$ 
remains unknown.
\begin{remark}
The expression on the right hand side of (\ref{T.n22}) is known (see \cite[p.168]{Ber96})
to be the Laplace transform of the random variable 
$$
T_{\max}:=\sup\{t>0:Y_{\alpha,\beta}(t)-t=\max_{u\geq0}(Y_{\alpha,\beta}(u)-u)\}.
$$
Let $f_{\max}$ denote the density function of this random variable. Then from (\ref{T.n22}) one can obtain the equality
$$
1-F_{\alpha,\beta}(x)=Cx^{\alpha(1-\rho)/(\alpha-1)}f_{\max}(x^{\alpha/(\alpha-1)}),\ x>0.
$$
Having this relation one can get the explicit form of $f_{\max}$ in the case of Brownian motion ($\alpha=2,\beta=0$) and in the case of spectrally
positive Levy processes ($\alpha\in(1,2),\beta=1$).
\hfill$\diamond$
\end{remark}

We now turn our attention to the moments of $\tau^{(a)}$.

It was shown by Gut \cite{Gut74} that the condition $\mathbf{E}(\max\{0,X\})^r<\infty$ for some $r>0$
is necessary and sufficient for the finiteness of $\mathbf{E}\bigl(\tau^{(a)}\bigr)^r$.
Therefore, the condition $X\in\mathcal{D}(\alpha,\beta)$ yields the finiteness of 
$\mathbf{E}\bigl(\tau^{(a)}\bigr)^r$ for all $r<\alpha$.

From the bound 
$$
\mathbf{P}(\tau^{(a)}>n)\leq\mathbf{P}(\tau^{(0)}>n)\text{ for all }n\geq0
$$ 
and (\ref{fixed}), using dominated convergence, we infer that
\begin{equation}
\label{small-r}
\lim_{a\to0}\mathbf{E}\bigl(\tau^{(a)}\bigr)^r=\mathbf{E}\bigl(\tau^{(0)}\bigr)^r<\infty
\end{equation}
for all $r\in(0,1-\rho)$.
Furthermore, it easy follows from Theorem~\ref{T.normal} and (\ref{tau0}) that
$$
\lim_{a\to0}\mathbf{E}\bigl(\tau^{(a)}\bigr)^r=\infty\text{ for all }r>1-\rho.
$$
Theorem~\ref{T.normal} allows us to determine the rate of growth as $a\to0$ of  $\mathbf{E}\bigl(\tau^{(a)}\bigr)^r$
for $r\in(1-\rho,\alpha)$. 
\begin{theorem}\label{T.exp}
Suppose $X\in\mathcal{D}(\alpha,\beta)$.
Then, for every $r\in(1-\rho,\alpha)$ there exists a function $L_r$ slowly varying at infinity such that
\begin{equation}
\label{T.n2}
\mathbf{E}\bigl(\tau^{(a)}\bigr)^r=L_r(1/a)a^{-\alpha(r+\rho-1)/(\alpha-1)}.
\end{equation}
\end{theorem}
This is already known in some particular cases, we now want to mention.

First of all we note that if the second moment of $X$ is finite, then, 
applying dominated convergence, one can show that 
$\mathbf{E}S^{(a)}_{\tau^{(a)}}\to\mathbf{E}S^{(0)}_{\tau^{(0)}}$ as $a\to0$. 
Thus, using the Wald identity and the well-known equality (see \cite[Proposition 18.5]{Sp01})
$$
-\mathbf{E}S^{(0)}_{\tau^{(0)}}=\frac{(\mathbf{E}X^2)^{1/2}}{\sqrt{2}}
\exp\Bigl\{\sum_{k=1}^\infty k^{-1}\Bigl(\mathbf{P}(S_k^{(0)}\geq0)-1/2\Bigr)\Bigr\},
$$ 
we obtain, as $a\to0$, 
\begin{equation}
\label{exp}
\mathbf{E}\tau^{(a)}\sim\frac{-\mathbf{E}S^{(0)}_{\tau^{(0)}}}{a}
=\frac{(\mathbf{E}X^2)^{1/2}}{a\sqrt{2}}\exp\Bigl\{\sum_{k=1}^\infty k^{-1}\Bigl(\mathbf{P}(S_k^{(0)}\geq0)-1/2\Bigr)\Bigr\}.
\end{equation}

Furthermore, the asymptotic behaviour of $\mathbf{E}\tau^{(a)}$ in the case of a non-Gaussian stable limit, that is, $\alpha<2$,
was recently studied by Lotov \cite{Lot06}. He has proved that 
$$
\mathbf{E}\tau^{(a)}=a^{-\alpha\rho/(\alpha-1)+o(1)}
\quad\text{as }a\downarrow0
$$
in this case.
Moreover, he has shown that (\ref{T.n2}) with $r=1$ holds under the additional condition
\begin{equation*}
\sum_{k=1}^\infty\frac{1}{k}\sup_{x\in\mathbb{R}}\big|\mathbf{P}(S^{(0)}_k>c_k x)-\mathbf{P}(Y_{\alpha,\beta}>x)\big|<\infty.
\end{equation*}

Having expressions for the expectation $\mathbf{E}\tau^{(a)}$ one can 
describe the asymptotic behaviour of some further characteristics of 
the random walk $\{S_n^{(a)},n\geq0\}$.
First, from the Wald identity and Theorem~\ref{T.exp} we obtain the equality
$$
\mathbf{E}S^{(a)}_{\tau^{(a)}}=-a\mathbf{E}\tau^{(a)}=-L_1(1/a)a^{1-\alpha\rho/(\alpha-1)}.
$$ 
Second, it is well known that the stopping time $\tau^{(a)}_+:=\min\{k\geq1:S^{(a)}_k\geq0\}$ 
is infinite with positive probability and $\mathbf{P}(\tau^{(a)}_+=\infty)=1/\mathbf{E}\tau^{(a)}$.
Then, using Theorem~\ref{T.exp} once again, we get
$$
\mathbf{P}(\tau^{(a)}_+=\infty)=a^{\alpha\rho/(\alpha-1)}/L_1(1/a).
$$

In conclusion of this subsection we note that our assumption that the distributions of $X^{(a)}$ and $X-a$ are equal can be weakened.
First of all we note, that if $X^{(a)}$ satisfies the conditions
$$
\mathbf{E}X^{(a)}=-a\quad\text{and}\quad\lim_{a\to0}\mathbf{E}\left(X^{(a)}\right)^2=\sigma^2\in(0,\infty),
$$
then the results of the present subsection are still hold.
Moreover, in the case of infinite second moment, the results of the present subsection remain valid if $X^{(a)}=X-a+Y^{(a)}$ in distribution,
where $X\in\mathcal{D}(\alpha,\beta)$ for some $\alpha\in(1,2)$ and $Y^{(a)}$ is such that
$$
\mathbf{E}Y^{(a)}=0,\quad Y^{(a)}\to0\text{ in law and }
\sup_{a\in[0,a_0]}\mathbf{E}\left|Y^{(a)}\right|^{\alpha+\delta}<\infty\text{ for some }\delta>0.
$$
We did not use these generalisations in the statements of our theorems because of results in the next subsection, where we need the assumption
$X^{(a)}=X-a$ in law.
\subsection{Results on large deviations}
If $na/c_n\to\infty$, then Theorem~\ref{T.normal} says only that
$$
\mathbf{P}(\tau^{(a)}>n)=o\bigl(\mathbf{P}(\tau^{(0)}>n)\bigr)\quad\text{as }a\to0.
$$
Our next purpose is to refine this relation and to find the rate of divergence
of $\mathbf{P}(\tau^{(a)}>n)$ in the mentioned above domain of {\it large deviations} for $\tau^{(a)}$.
 To proceed in this situation one has to know the
asymptotic behaviour of $\mathbf{P}(S_n^{(a)}>0)$ for $na/c_n\to\infty$.
It follows from the definition of $S_n^{(a)}$ that 
$\mathbf{P}(S_n^{(a)}>0)=\mathbf{P}(S_n^{(0)}>na)$. Thus, the assumption
$na/c_n\to\infty$ means that we are in the domain  of large
deviations for $S_n^{(0)}$. Since the behaviour of large deviation probabilities 
depends crucially on whether the limit of $S_n^{(0)}/c_n$ is Gaussian
or strictly stable,i.e, $\alpha\in(1,2)$, we consider these two cases separately.

If $S_n^{(0)}$ belongs to the domain of attraction of a strictly stable law, 
then, as is well known, 
$$
\mathbf{P}(S_n^{(0)}\geq x_n)\sim n\mathbf{P}(X\geq x_n)
$$ 
for any sequence $x_n$ satisfying $x_n/c_n\to\infty$. This relation allows 
one to obtain the following result.
\begin{theorem}\label{LD_stable}
Suppose $X\in\mathcal{D}(\alpha,\beta)$ for some $1<\alpha<2$ and $\beta>-1$.\\
If $n=n(a)$ is such that $na/c_n\to\infty$, then
\begin{equation}
\label{T.n3}
\mathbf{P}(\tau^{(a)}>n)\sim\mathbf{E}\tau^{(a)}\mathbf{P}(X\geq na)\quad\text{as }a\to0.
\end{equation}
\end{theorem}
The right hand side of (\ref{T.n3}) coincides with that of (\ref{Don1}). 
Roughly speaking, if $n$ is very large, then the asymptotic behaviour of 
$\mathbf{P}(\tau^{(a)}>n)$ for $a\to0$ is as in the case of the fixed negative
drift. But there is one crucial difference between fixed and asymptotically small drift:
The expectation $\mathbf{E}\tau^{(a)}$ grows unbounded if $a\to0$, and is
a constant when the drift is fixed. Therefore, (\ref{T.n3}) would be
useless without Theorem~\ref{T.exp}.

We turn our attention to the case when $\sigma^2:=\mathbf{E}X^2$ is finite. 
Here we shall assume, without loss of generality, that $\sigma^2=1$.
Under this condition we have $c_n=\sqrt{n}$. Then the condition 
$an/c_n\to\infty$ reads as $na^2\to\infty$. In this case of finite variance
the asymptotic behaviour of $\mathbf{P}(S_n^{(0)}>x_n)$ depends not only
on the tail behaviour of $X$, but also on the rate of the growth of $x_n$. 
If $x_n$ grows not very fast
($x_n=o(r_1(n))$ for some $r_1(n)$ depending on the distribution of $X$), 
then one has an asymptotic expression for $\mathbf{P}(S_n^{(0)}>x_n)$
in terms of the so-called Cram\'er series (for the definition of the
Cram\'er series see, for example, \cite[Chapter VIII]{Pet75}). 
For this type of large deviations we have the following result.
\begin{theorem}\label{LD_var}
Assume that $\mathbf{E}X^2=1$, $n=n(a)$ is such that $na^2\to\infty$, and that
\begin{equation}
\label{cond1}
\mathbf{P}(S_j^{(0)}\geq ja)\sim \overline{\Phi}(\sqrt{j}a)
\exp\{ja^3\lambda_m(a)\}
\text{ uniformly in }j\in[a^{-2},n],
\end{equation}
where $\lambda_m(u)$ is the partial sum in the Cram\'er series containing the first $m$ terms
and $\overline{\Phi}(x):=\int_x^\infty\frac{1}{\sqrt{2\pi}}e^{-u^2/2}du$.
Then
\begin{equation}
\label{T.n4}
\mathbf{P}(\tau^{(a)}>n)\sim 2\mathbf{E}\tau^{(a)}\frac{1}{n}\overline{\Phi}(\sqrt{n}a)
\exp\{na^3\lambda_m(a)\}.
\end{equation}
\end{theorem}
Condition (\ref{cond1}) has one essential disadvantage: it involves the whole
sequence $\{S_k^{(0)},k\geq0\}$. We now list some restrictions on the distribution of $X$, 
which imply the validity of (\ref{cond1}). 

Nagaev S.V. \cite{Nag65} has
proved that the condition $\mathbf{E}|X|^k<\infty$ with some $k>2$ implies that the relation
\begin{equation}
\label{NormApp}
\mathbf{P}(S_n^{(0)}\geq x)\sim\overline{\Phi}(x/\sqrt{n})\quad\text{as }n\to\infty
\end{equation} 
holds uniformly in $x\leq \sqrt{(k/2-1)n\log n}$. Thus, the existence of $\mathbf{E}|X|^k$ for some $k>2$
yields (\ref{cond1}) with $m=0$ for all $n$ satisfying 
$$
n\leq\Bigl(\frac{k}{2}-1\Bigr)a^{-2}\log a^{-2}.
$$
Furthermore, it has been proved by Nagaev A.V. \cite{Nag69} and by Rozovskii \cite{Roz89}
that if $\mathbf{P}(X>x)$ is regularly varying at infinity with
index $p<-2$, then, under some additional restrictions on the left tail, 
\begin{equation}
\label{NagRoz}
\mathbf{P}(S_n^{(0)}\geq x)\sim\overline{\Phi}(x/\sqrt{n})+n\mathbf{P}(X>x+\sqrt{n})
\end{equation}
uniformly on $x>0$. 
Thus, (\ref{NormApp}) holds for
all $x\leq C\sqrt{n\log n}$ for any $C<(p-2)^{1/2}$. Consequently, (\ref{cond1})
with $m=0$ holds for
$$
n\leq Ca^{-2}\log a^{-2},\ C<(p-2)^{1/2}.
$$

Osipov \cite{Os72} has found necessary and sufficient conditions, under which 
the relation
$$
\mathbf{P}(S_n^{(0)}\geq x)\sim \overline{\Phi}(x/\sqrt{n})
\exp\Bigl\{\frac{x^3}{n^2}\lambda_{[1/(1-\gamma)]}\Bigl(\frac{x}{n}\Bigr)\Bigr\}
$$
holds uniformly in $0\leq x\leq n^{\gamma},\ 1/2<\gamma<1$, where $[t]$ denotes the integer part
of $t$. If these conditions are fulfilled, then, obviously, (\ref{cond1}) holds with 
$m=[1/(1-\gamma)]$ for all $n\leq a^{1/(1-\gamma)}$.

It is well-known that
if $X$ satisfies the Cram\'er condition ($\mathbf{E}e^{h|X|}<\infty$ for some $h>0$), 
then (\ref{cond1}) holds with $m=\infty$ and for all $n$ satisfying $na^2\to\infty$. 
Thus, Theorems~\ref{T.normal} and \ref{LD_var} describe the behaviour of 
$\mathbf{P}(\tau^{(a)}>n)$ for any choice of $n=n(a)$ and any random walk satisfying 
the Cram\'er condition. 

It is easy to see that the statement of Theorem~\ref{LD_var} can be rewritten as follows:
If (\ref{cond1}) holds, then
\begin{equation*}
\mathbf{P}(\tau^{(a)}>n)\sim \frac{2}{\sqrt{2\pi}}a^{-1}\mathbf{E}\tau^{(a)}
n^{-3/2}e^{-n\xi(a)},
\end{equation*}
where 
\begin{equation}
\label{xi-def}
\xi(a):=\frac{a^2}{2}-a^3\lambda_m(a).
\end{equation}
Furthermore, in the proof of Theorem~\ref{LD_var} we shall see that
\begin{equation*}
\frac{\mathbf{E}[e^{\xi(a)\tau^{(a)}},\tau^{(a)}\leq n]-1}{e^{\xi(a)}-1}\sim 2\mathbf{E}\tau^{(a)}.
\end{equation*}
Thus, 
\begin{equation*}
\mathbf{P}(\tau^{(a)}>n)\sim \frac{1}{a\sqrt{2\pi}}
\frac{\mathbf{E}[e^{\xi(a)\tau^{(a)}},\tau^{(a)}\leq n]-1}{e^{\xi(a)}-1}
n^{-3/2}e^{-n\xi(a)},
\end{equation*}
which is rather close to relation (\ref{Don2}). If, additionally, $X$ satisfies
the Cram\'er condition, implying (\ref{cond1}) with $m=\infty$, then one can replace
the truncated expectation $\mathbf{E}[e^{\xi(a)\tau^{(a)}},\tau^{(a)}\leq n]$ by
$\mathbf{E}[e^{\xi(a)\tau^{(a)}}]$:
\begin{equation}
\label{cor}
\mathbf{P}(\tau^{(a)}>n)\sim \frac{1}{a\sqrt{2\pi}}
\frac{\mathbf{E}[e^{\xi(a)\tau^{(a)}}]-1}{e^{\xi(a)}-1}
n^{-3/2}e^{-n\xi(a)}.
\end{equation}
It follows from the definition of the Cram\'er series that $\xi(a)$, defined in (\ref{xi-def}), is the unique positive 
solution to the equation $\frac{d}{dh}\mathbf{E}e^{hX^{(a)}}=0$. Therefore, (\ref{cor})
is an analog of (\ref{Don2}) for random walks with vanishing drift.

Another type of large deviation behaviour appears in the case when $x_n$ grows fast, i.e.,
$x_n\gg r_2(n)$ and the tail of $X$ varies in an appropriate way.
(Recall that $a_n\gg b_n$ means that $\frac{a_n}{b_n}\to\infty$.) Here, as in the
case of non-gaussian stable limit, one has $\mathbf{P}(S_n^{(0)}\geq x_n)\sim n\mathbf{P}(X\geq x_n)$.
We consider only the case when the tail of $X$ is regularly varying.
\begin{theorem}
\label{tail}
Assume that $\mathbf{P}(X\geq x)$ is regularly varying at infinity with index $p<-2$ and
\begin{equation}
\label{Roz_cond}
\int_{|x|>y}x^2\mathbf{P}(X\in dx)=o\Bigl(\frac{1}{\log y}\Bigr)\text{ as }y\to\infty.
\end{equation}
Then, as $a\to0$,
$$
\mathbf{P}(\tau^{(a)}>n)\sim\mathbf{E}\tau^{(a)}\mathbf{P}(X\geq na)
$$
for any $n=n(a)$ satisfying the inequality $n(a)\geq Ca^{-2}\log a^{-2}$ with some $C>(p-2)^{1/2}$.
\end{theorem}
After Theorem~\ref{LD_var} we have mentioned that, in the case of regularly varying tails,
(\ref{cond1}) holds for all $n\leq Ca^{-2}\log a^{-2}$, $C<(p-2)^{1/2}$. Therefore, the
behaviour of $\mathbf{P}(\tau^{(a)}>n)$ remains unclear only for $n$ satisfying
$(na^2/\log a^{-2})\to (p-2)^{1/2}$. We conjecture that if the conditions of Theorem~\ref{tail}
hold, then, in agreement with (\ref{NagRoz}),
$$
\mathbf{P}(\tau^{(a)}>n)\sim2\mathbf{E}\tau^{(a)}\frac{1}{n}\overline{\Phi}(\sqrt{n}a)+
\mathbf{E}\tau^{(a)}\mathbf{P}(X\geq\sqrt{n}+na)
$$
for all $n$ satisfying $na^2\to\infty$.

The remaining part of the paper is organised as follows. In the next section we derive an upper bound for
the probability $\mathbf{P}(\tau^{(a)}>n)$, which is crucial for the proof of Theorem~\ref{T.exp}. 
This proof will be given in Section~\ref{s_T.exp}. Section~\ref{s_T.normal} is devoted to the proof of
Theorem~\ref{T.normal}. Finally, Theorems \ref{LD_stable}, \ref{LD_var} and \ref{tail} will be proved in the
last Section.
\section{Upper bounds for the tail of $\tau^{(a)}$}
It follows from (\ref{GenF}) that in order to obtain upper 
bounds for $\mathbf{P}(\tau^{(a)}>n)$ one needs inequalities 
for $\mathbf{P}(S_n^{(a)}\geq 0)=\mathbf{P}(S_n^{(0)}\geq na)$.
In the following lemma we adapt one of the well-known Fuk-Nagaev
inequalities for our purposes.
\begin{lemma}
\label{Fuk-Nagaev}
Assume that $X\in\mathcal{D}(\alpha,\beta)$.
Then there exists a constant $C$ such that the inequality
\begin{equation*}
\mathbf{P}(S_n^{(0)}\geq x)\leq n\mathbf{P}(X\geq x/3)+C\Bigl(\frac{nV(x)}{x^2}\Bigr)^2
\end{equation*}
holds for all $x>0$ and $n\geq1$.
\end{lemma}
\begin{proof}
Applying Theorem ~1.2 of \cite{Nag79} with $t=2$, we have
\begin{equation}
\label{FN1}
\mathbf{P}(S_n^{(0)}\geq x)\leq n\mathbf{P}(X\geq y)+e^{x/y}
\Bigl(\frac{nV(y)}{xy}\Bigr)^{x/y+nV(y)/y^2-n\mu(y)/y},
\end{equation}
where $\mu(y):=\mathbf{E}[X,\,|X|\leq y]$.

Since $\mathbf{E}X=0$,
\begin{align*}
|\mu(y)|=\Big|\int_{|x|>y}x\mathbf{P}(X\in dx)\Big|&\leq \int_{x>y}x\mathbf{P}(|X|\in dx)\\
&=y\mathbf{P}(|X|>y)+\int_y^\infty\mathbf{P}(|X|>x)dx.
\end{align*}
It is well-known that the assumption $X\in\mathcal{D}(\alpha,\beta)$ yields
$$
\lim_{x\to\infty}\frac{x^2\mathbf{P}(|X|>x)}{V(x)}=\frac{2-\alpha}{\alpha}.
$$
Therefore, as $y\to\infty$,
\begin{align*}
|\mu(y)|&\leq \left(\frac{2-\alpha}{\alpha}+o(1)\right)\left(\frac{V(y)}{y}+\int_y^\infty\frac{V(x)}{x^2}dx\right)\\
&=\left(\frac{2-\alpha}{\alpha-1}+o(1)\right)\frac{V(y)}{y}.
\end{align*}
In the last step we used the relation
$$
\int_y^\infty\frac{V(x)}{x^2}dx\sim\frac{1}{\alpha-1}\frac{V(y)}{y}\quad\text{as }y\to\infty,
$$
which follows from the fact that $V(x)$ is regularly varying with index $2-\alpha$.
As a result we have the bound
\begin{equation}
\label{bound}
\frac{V(y)}{y^2}-\frac{\mu(y)}{y}\geq\left(\frac{2\alpha-3}{\alpha-1}+o(1)\right)\frac{V(y)}{y^2}.
\end{equation}
It follows from definition (\ref{Defa}) of the sequence $\{c_n\}$ that $V(c_n)/c_n^2\sim n^{-1}$ as $n\to\infty$.
Consequently, there exists a constant $C(\alpha)$ such that
\begin{equation*}
\frac{V(y)}{y^2}-\frac{\mu(y)}{y}\geq -\frac{1}{n}
\end{equation*}
for all $y>C(\alpha)c_n$. From this bound and (\ref{FN1}) with $y=x/3$ we get
\begin{equation*}
\mathbf{P}(S_n^{(0)}\geq x)\leq n\mathbf{P}(X\geq x/3)+27e^{3}
\Bigl(\frac{nV(y)}{x^2}\Bigr)^{2},\quad x\geq 3C(\alpha)c_n.
\end{equation*}
This inequality, together with monotonicity of $V$, implies that
the desired result holds for $x>C(\alpha)c_n$. Noting that 
\begin{equation*}
\min_{n\geq1}\inf_{x\leq 3C(\alpha)c_n}\frac{nV(x)}{x^2}>0,
\end{equation*}
we complete the proof of the lemma.
\end{proof}

In order to `translate' bounds for $\mathbf{P}(S_n^{(0)}>na)$
into bounds for $\mathbf{P}(\tau^{(a)}>n)$ we shall use the recurrent relation
\begin{equation}
\label{Rec}
n\mathbf{P}(\tau^{(a)}>n)=\sum_{j=0}^{n-1}\mathbf{P}(\tau^{(a)}>j)\mathbf{P}\left(S_{n-j}^{(0)}>(n-j)a\right),
\end{equation}
which can be obtained by differentiating (\ref{GenF}).
\begin{proposition}
\label{UpBound}
The inequality 
\begin{equation*}
\mathbf{P}(\tau^{(a)}>n)\leq C\mathbf{E}\tau^{(a)}\frac{V(na)}{(na)^2}
\end{equation*}
is valid for all $a>0$ and all $n\geq n_a:=\min\{n\geq1:an>c_n\}$.
\end{proposition}
\begin{proof}
Using Lemma~\ref{Fuk-Nagaev}, we have
\begin{align}
\label{P1}
\nonumber
&\sum_{0\leq j<n/2}\mathbf{P}(\tau^{(a)}>j)\mathbf{P}\left(S_{n-j}^{(0)}>(n-j)a \right)\\
\nonumber
&\hspace{0.5cm}
\leq\sum_{0\leq j<n/2}\mathbf{P}(\tau^{(a)}>j)\Bigl((n-j)\mathbf{P}\left(X\geq (n-j)a/3\right)+
C\Bigl(\frac{(n-j)V((n-j)a)}{((n-j)a)^2}\Bigr)^{2}\Bigr)\\
\nonumber
&\hspace{0.5cm}\leq\Bigl(n\mathbf{P}(X\geq na/6)+C\Bigl(\frac{nV(na)}{(na)^2}\Bigr)^{2}\Bigr)
\sum_{0\leq j<n/2}\mathbf{P}(\tau^{(a)}>j)\\
\nonumber
&\hspace{0.5cm}\leq\mathbf{E}\tau^{(a)}\Bigl(n\mathbf{P}(X\geq na/6)+C\Bigl(\frac{nV(na)}{(na)^2}\Bigr)^{2}\Bigr)\\
&\hspace{0.5cm}\leq n\mathbf{E}\tau^{(a)}\Bigl(\mathbf{P}(X\geq na/6)+C\frac{V(na)}{(na)^2}\Bigr).
\end{align}
In the last step we used definition (\ref{Defa}) of $c_n$ and the bound $an\geq c_n$, which follows from the assumption $n\geq n_a$.

Further, using the Markov inequality, we get
\begin{equation}
\label{P2}
\sum_{n/2\leq j\leq n-1}\mathbf{P}(\tau^{(a)}>j)\mathbf{P}\left(S_{n-j}^{(0)}>(n-j)a\right)\leq\frac{2\mathbf{E}\tau^{(a)}}{n}
\sum_{k=1}^n\mathbf{P}\left(S_k^{(0)}\geq ka\right).
\end{equation}
Applying Lemma~\ref{Fuk-Nagaev}, we obtain
\begin{align}
\label{P3}
\nonumber
\sum_{k=1}^n\mathbf{P}(S_k^{(0)}\geq ka)&\leq n_a+\sum_{k=n_a}^n\mathbf{P}(S_k^{(0)}\geq ka)\\
&\leq n_a+\sum_{k=n_a}^nk\mathbf{P}(X\geq ka/3)+C\sum_{k=n_a}^n\frac{V^2(ka)}{k^2a^4}.
\end{align}
Since $V(x)$ is regularly varying with index $2-\alpha$,
\begin{align}
\label{P4}
\nonumber
\sum_{k=n_a}^n\frac{V^2(ka)}{k^2a^4}&\leq Ca^{-2}\sum_{k=n_a}^n\frac{V^2(ka)}{(ka)^2}\leq Ca^{-3}\int_{an_a}^{an}\frac{V^2(x)}{x^2}dx\\
&\leq Ca^{-3}V(an)\int_{an_a}^{an}\frac{V(x)}{x^2}dx\leq Ca^{-3}V(an)\frac{V(an_a)}{an_a}.
\end{align}
From the definitions of $c_n$ and $n_a$ we infer that 
\begin{equation}
\label{P5}
V(an_a)\sim V(c_{n_a})\sim\frac{c_{n_a}^2}{n_a}\sim a^2n_a.
\end{equation}
Applying this relation to the last line in array (\ref{P4}), we obtain the bound
\begin{equation}
\label{P6}
\sum_{k=n_a}^n\frac{V^2(ka)}{k^2a^4}\leq Ca^{-2}V(an).
\end{equation}
Furthermore,
\begin{align}
\label{P7}
\nonumber
\sum_{k=n_a}^nk\mathbf{P}(X\geq ka/3)&\leq Ca^{-2}\int_{an_a}^{an}x\mathbf{P}(X>x/3)dx\\
\nonumber
&\leq a^{-2}\int_0^{na}x\mathbf{P}(|X|>x)dx\\
&=\frac{a^{-2}}{2}\Bigl(V(an)+(an)^2\mathbf{P}(|X|>an)\Bigr),
\end{align}
where in the last step we used integration by parts. Combining (\ref{P3}), (\ref{P6}) and (\ref{P7}), we have
\begin{equation}
\label{P8}
\sum_{k=1}^n\mathbf{P}(S_k^{(0)}\geq ka)\leq Cn_a+ Ca^{-2}V(an)+n^2\mathbf{P}(|X|>an).
\end{equation}
It is easy to see that (\ref{P5}) yields $n_a\sim a^{-2}V(an_a)$. From this relation and monotonicity of $V(x)$ we conclude
that $n_a\leq Ca^{-2}V(an)$ for all $n\geq n_a$. Applying this bound to (\ref{P8}), we get
\begin{equation}
\label{P9}
\sum_{k=1}^n\mathbf{P}(S_k^{(0)}\geq ka)\leq Ca^{-2}V(an)+n^2\mathbf{P}(|X|>an).
\end{equation}
Combining (\ref{P1}), (\ref{P2}) and (\ref{P9}), we arrive at the inequality
\begin{equation}
\label{P10}
\sum_{j=0}^{n}\mathbf{P}(\tau^{(a)}>j)\mathbf{P}\left(S_{n-j}^{(0)}>(n-j)a\right)\leq
C n\mathbf{E}\tau^{(a)}\Bigl(\mathbf{P}(|X|\geq na/6)+\frac{V(na)}{(na)^2}\Bigr).
\end{equation}
It is easy to see that
\begin{align*}
\mathbf{P}(|X_1|>x)&=\sum_{j=0}^\infty\mathbf{P}\left(|X_1|\in(2^jx,2^{j+1}x]\right)\leq
\sum_{j=0}^\infty\frac{V(2^{j+1}x)}{2^{2j}x^2}\\
&\leq\frac{V(x)}{x^2}4C(\gamma)\sum_{j=1}^\infty 2^{-(\alpha-\gamma)j}.
\end{align*}
Here we used the inequality
\begin{equation}
\label{eq:L2a}
\frac{V(y)}{V(x)}\leq C(\gamma)\Bigl(\frac{y}{x}\Bigr)^{2-\alpha+\gamma},\quad y\geq x,
\end{equation}
which follows from the Karamata representation, see \cite[Theorem 1.2]{Sen76}, recall that
$V(x)$ is regularly varying with index $2-\alpha$. Choosing $\gamma<\alpha$, we get
\begin{equation*}
\mathbf{P}(|X_1|>x)\leq C\frac{V(x)}{x^2}.
\end{equation*}
Therefore, the right hand side in (\ref{P10}) is bounded by $Cn\mathbf{E}\tau^{(a)}\frac{V(na)}{(na)^2}$.
Thus, the statement of the proposition follows from (\ref{Rec}).
\end{proof}
\section{Proof of Theorem~\ref{T.normal}}\label{s_T.normal}
From the definition of the first ladder epoch $\tau^{(a)}$ we get
\begin{align}
\label{Rep}
\nonumber
\mathbf{P}(\tau^{(a)}>n)&=\mathbf{P}\Bigl(\min_{1\leq k\leq n}(S_k^{(0)}-ka)>0\Bigr)\\
\nonumber
&=\mathbf{P}\Bigl(\min_{1\leq k\leq n}S_k^{(0)}>0\Bigr)
\mathbf{P}\Bigl(\min_{1\leq k\leq n}(S_k^{(0)}-ka)>0\Big|\min_{1\leq k\leq n}S_k^{(0)}>0\Bigr)\\
&=\mathbf{P}(\tau^{(0)}>n)
\mathbf{P}\Bigl(\min_{1\leq k\leq n}\Bigl(\frac{S_k^{(0)}}{c_n}-\frac{k}{n}\frac{an}{c_n}\Bigr)>0\Big|\min_{1\leq k\leq n}S_k^{(0)}>0\Bigr).
\end{align}
Doney \cite{Don85} has shown that $\{S_{[tc_n]}^{(0)}/c_n,\ t\in[0,1]|\min_{1\leq k\leq n}S_k^{(0)}>0\}$ converges
weakly, as $n\to\infty$, to the Levy meander $\{M_{\alpha,\beta}(t),\ t\in[0,1]\}$. This yields
\begin{align}
\label{Invar}
\nonumber
\lim_{n\to\infty}
\mathbf{P}\Bigl(\min_{1\leq k\leq n}\Bigl(\frac{S_k^{(0)}}{c_n}-\frac{k}{n}\frac{an}{c_n}\Bigr)>0\Big|\min_{1\leq k\leq n}S_k^{(0)}>0\Bigr)
&=\mathbf{P}\Bigl(\min_{0\leq t\leq1}(M_{\alpha,\beta}(t)-ut)>0\Bigr)\\
&=:1-F_{\alpha,\beta}(u).
\end{align}
It is obvious that $F_{\alpha,\beta}(u)$ is monotonously increasing and $\lim_{u\to\infty}F_{\alpha,\beta}(u)=1$. 

It is known that the corresponding meander $M_{\alpha,\beta}$ can be defined by
$$
\{M_{\alpha,\beta}(t),\ t\in[0,1]\}=\lim_{\varepsilon\to0}
\{Y_{\alpha,\beta}(t),\ t\in[0,1]|\inf_{0\leq t\leq1}Y_{\alpha,\beta}(t)>0,Y_{\alpha,\beta}(0)=\varepsilon\}.
$$ 
Therefore,
$$
1-F_{\alpha,\beta}(u)=\lim_{\varepsilon\to0}\frac{\mathbf{P}\Bigl(\inf_{0\leq t\leq1}(Y_{\alpha,\beta}(t)-ut)>0|Y_{\alpha,\beta}(0)=\varepsilon\Bigr)}
{\mathbf{P}\Bigl(\inf_{0\leq t\leq1}Y_{\alpha,\beta}(t)>0|Y_{\alpha,\beta}(0)=\varepsilon\Bigr)}.
$$
Define $H_{\alpha,\beta}^{(u)}(z):=\min\{t:Y_{\alpha,\beta}(t)-ut\leq z|Y_{\alpha,\beta}(0)=0\}$. Then
$$
1-F_{\alpha,\beta}(u)=
\lim_{\varepsilon\to0}\frac{\mathbf{P}(H_{\alpha,\beta}^{(u)}(-\varepsilon)>1)}{\mathbf{P}(H_{\alpha,\beta}^{(0)}(-\varepsilon)>1)}.
$$
In the case of the Brownian motion, that is, $\alpha=2,\beta=0$,  one can calculate the limit explicitly. Indeed, it is known that 
$H_{2,0}^{(u)}(-\varepsilon)$ has the density
$$
\frac{\varepsilon}{\sqrt{2\pi}t^{3/2}}\exp\Bigl\{-\frac{(ut-\varepsilon)^2}{2t}\Bigr\},\ t>0.
$$
Thus, as $\varepsilon\to0$,
$$
\mathbf{P}(H_{2,0}^{(0)}(-\varepsilon)>1)=\frac{\varepsilon}{\sqrt{2\pi}}\int_1^\infty t^{-3/2}e^{-\varepsilon^2/2t}dt
\sim\frac{2\varepsilon}{\sqrt{2\pi}},
$$
and, consequently,
\begin{align*}
\lim_{\varepsilon\to0}\frac{\mathbf{P}(H_{\alpha,\beta}^{(u)}(-\varepsilon)>1)}{\mathbf{P}(H_{\alpha,\beta}^{(0)}(-\varepsilon)>1)}
&=\lim_{\varepsilon\to0}\int_1^\infty\frac{1}{2t^{3/2}}\exp\Bigl\{-\frac{(ut-\varepsilon)^2}{2t}\Bigr\}dt\\
&=\int_1^\infty\frac{1}{2t^{3/2}}e^{-u^2t/2}dt=u\int_u^\infty v^{-2}e^{-v^2/2}dv.
\end{align*}
As a result we have 
\begin{equation}
\label{F6}
1-F_{2,0}(u)=u\int_u^\infty v^{-2}e^{-v^2/2}dv.
\end{equation}

This equality can be generalised to stable Levy processes without negative jumps, i.e., 
$\{\alpha\in(1,2),\beta=1\}$ or $\{\alpha=2,\beta=0\}$.
Indeed, using Kendall's equality (see \cite{Keil73}) and the scaling property of stable processes, 
we see that $H_{\alpha,1}^{(u)}(-\varepsilon)$ has the density
$$
u\mapsto\frac{\varepsilon}{t^{1+1/\alpha}}g_{\alpha,1}\Bigl(\frac{-\varepsilon+ut}{t^{1/\alpha}}\Bigr).
$$
Then, analogously to the case of the Brownian motion,
$$
1-F_{\alpha,1}(u)=\frac{u^{1/(\alpha-1)}}{(\alpha-1)g_{\alpha,1}(0)}\int_u^\infty v^{-\alpha/(\alpha-1)}g_{\alpha,1}(v)dv.
$$

Unfortunately we can not give an explicit expression for $1-F_{\alpha,\beta}$ for a process with positive jumps.
But we can describe this function via Laplace transform of $x^{\rho-1}(1-F_{\alpha,\beta}(x^{1-1/\alpha}))$.

In order to prove (\ref{T.n22}) we show that $1-F_{\alpha,\beta}$ satisfies a certain integral equation.
Dividing both parts of (\ref{Rec}) by $n\mathbf{P}(\tau^{(0)}>n)$, we have
\begin{equation}
\label{n1}
\frac{\mathbf{P}(\tau^{(a)}>n)}{\mathbf{P}(\tau^{(0)}>n)}=
\sum_{j=0}^{n-1}\frac{\mathbf{P}(\tau^{(a)}>j)}{\mathbf{P}(\tau^{(0)}>j)}
\frac{\mathbf{P}(\tau^{(0)}>j)}{\mathbf{P}(\tau^{(0)}>n)}\mathbf{P}\left(S^{(0)}_{n-j}\geq a(n-j)\right)\frac{1}{n}.
\end{equation}
Fix any $\varepsilon\in(0,1/2)$. We first note that
\begin{align}
\label{n2}
\nonumber
\sum_{0\leq j\leq \varepsilon n}\frac{\mathbf{P}(\tau^{(a)}>j)}{\mathbf{P}(\tau^{(0)}>j)}
\frac{\mathbf{P}(\tau^{(0)}>j)}{\mathbf{P}(\tau^{(0)}>n)}&\mathbf{P}\left(S^{(0)}_{n-j}\geq a(n-j)\right)\frac{1}{n}\\
&\leq\frac{\sum\limits_{0\leq j\leq \varepsilon n}\mathbf{P}(\tau^{(0)}>j)}{n\mathbf{P}(\tau^{(0)}>n)}
\leq C\varepsilon^\rho
\end{align}
and
\begin{align}
\label{n3}
\nonumber
\sum_{(1-\varepsilon)n\leq j\leq n-1}\frac{\mathbf{P}(\tau^{(a)}>j)}{\mathbf{P}(\tau^{(0)}>j)}
\frac{\mathbf{P}(\tau^{(0)}>j)}{\mathbf{P}(\tau^{(0)}>n)}\mathbf{P}\left(S^{(0)}_{n-j}\geq a(n-j)\right)\frac{1}{n}\\
\leq\frac{\mathbf{P}(\tau^{(0)}>n/2)}{n\mathbf{P}(\tau^{(0)}>n)}\varepsilon n\leq C\varepsilon.
\end{align}
In both bounds we have used the fact that $\mathbf{P}(\tau^{(0)}>j)$ varies regularly with index $\rho-1$.

It remains to consider the middle part of the sum on the right hand side of (\ref{n1}). It is easy to see that
the condition $an/c_n\to u$ implies that 
$$
aj/c_j\to ut^{1-1/\alpha}\quad\text{as }a\to0,
$$ 
provided that $j\sim tn$. Then, in view of (\ref{Invar}),
for every $t\in(0,1)$ the following is valid: As $a\to0$,
\begin{align*}
f_a(t)&:=\frac{\mathbf{P}(\tau^{(a)}>[tn])}{\mathbf{P}(\tau^{(0)}>[tn])}
\frac{\mathbf{P}(\tau^{(0)}>[tn])}{\mathbf{P}(\tau^{(0)}>n)}\mathbf{P}\left(S^{(0)}_{n-[tn]}\geq a(n-[tn])\right)\\
&\to (1-F_{\alpha,\beta}(ut^{1-1/\alpha}))t^{\rho-1}\mathbf{P}\left(Y_{\alpha,\beta}(1)>u(1-t)^{1-1/\alpha}\right).
\end{align*}
Thus, by dominated convergence,
\begin{align*}
\lim_{a\to0}\sum_{\varepsilon n<j<(1-\varepsilon)n}\frac{\mathbf{P}(\tau^{(a)}>j)}{\mathbf{P}(\tau^{(0)}>j)}
\frac{\mathbf{P}(\tau^{(0)}>j)}{\mathbf{P}(\tau^{(0)}>n)}\mathbf{P}\left(S^{(0)}_{n-j}\geq a(n-j)\right)\frac{1}{n}\hspace{2cm}\\
=\int_\varepsilon^{1-\varepsilon}
(1-F_{\alpha,\beta}(ut^{1-1/\alpha}))t^{\rho-1}\mathbf{P}\left(Y_{\alpha,\beta}(1)>u(1-t)^{1-1/\alpha}\right)dt.
\end{align*}
Using now monotone convergence, we obtain
\begin{align}
\label{n4}
\nonumber
\lim_{\varepsilon\to0}\lim_{a\to0}\sum_{\varepsilon n<j<(1-\varepsilon)n}\frac{\mathbf{P}(\tau^{(a)}>j)}{\mathbf{P}(\tau^{(0)}>j)}
\frac{\mathbf{P}(\tau^{(0)}>j)}{\mathbf{P}(\tau^{(0)}>n)}\mathbf{P}\left(S^{(0)}_{n-j}\geq a(n-j)\right)\frac{1}{n}\hspace{1cm}\\
=\int_0^{1}
(1-F_{\alpha,\beta}(ut^{1-1/\alpha}))t^{\rho-1}\mathbf{P}\left(Y_{\alpha,\beta}(1)>u(1-t)^{1-1/\alpha}\right)dt.
\end{align}
Combining (\ref{n1}) -- (\ref{n4}), and taking into account (\ref{Invar}), we get
\begin{equation}
\label{n5}
1-F_{\alpha,\beta}(u)=\int_0^{1}
(1-F_{\alpha,\beta}(ut^{1-1/\alpha}))t^{\rho-1}\mathbf{P}\left(Y_{\alpha,\beta}(1)>u(1-t)^{1-1/\alpha}\right)dt.
\end{equation}

Setting
$$
G_{\alpha,\beta}(u):=1-F_{\alpha,\beta}(u^{1-1/\alpha})\text{ and }\xi_{\alpha,\beta}:=(Y_{\alpha,\beta}(1))^{\alpha/(\alpha-1)},
$$
we can rewrite (\ref{n5}) as follows
$$
G_{\alpha,\beta}(u)=\int_0^{1}
G_{\alpha,\beta}(ut)t^{\rho-1}\mathbf{P}\bigl(\xi_{\alpha,\beta}>u(1-t)\bigr)dt.
$$
Substituting $t=y/u$, we have
$$
G_{\alpha,\beta}(u)=u^{-\rho}\int_0^u
G_{\alpha,\beta}(y)y^{\rho-1}\mathbf{P}(\xi_{\alpha,\beta}>u-y)dy.
$$
Therefore, the function $Q_{\alpha,\beta}(u):=u^{\rho-1}G_{\alpha,\beta}(u)$ satisfies the equation
\begin{equation}
\label{n6}
uQ_{\alpha,\beta}(u)=\int_0^uQ_{\alpha,\beta}(y)\mathbf{P}(\xi_{\alpha,\beta}>u-y)dy.
\end{equation}
Let $q_{\alpha,\beta}(\lambda)$ denote the Laplace transform of the function $Q_{\alpha,\beta}$, i.e.,
$$
q_{\alpha,\beta}(\lambda)=\int_0^\infty e^{-\lambda u}Q_{\alpha,\beta}(u)du,\quad \lambda>0.
$$
Now(\ref{n6}) implies that
\begin{align*}
\frac{d}{d\lambda}q_{\alpha,\beta}(\lambda)&=-\int_0^\infty ue^{-\lambda u}Q_{\alpha,\beta}(u)du\\
&=-\int_0^\infty e^{-\lambda u}\int_0^uQ_{\alpha,\beta}(y)\mathbf{P}(\xi_{\alpha,\beta}>u-y)dy\\
&=-\int_0^\infty e^{-\lambda u}Q_{\alpha,\beta}(u)du\int_0^\infty e^{-\lambda z}\mathbf{P}(\xi_{\alpha,\beta}>z)dz\\
&=-q_{\alpha,\beta}(\lambda)\int_0^\infty e^{-\lambda z}\mathbf{P}(\xi_{\alpha,\beta}>z)dz.
\end{align*}
Solving this differential equation, we see that
\begin{align*}
q_{\alpha,\beta}(\lambda)&=q_{\alpha,\beta}(\lambda_0)\exp\Bigl\{
-\int_{\lambda_0}^\lambda\int_0^\infty e^{-\lambda z}\mathbf{P}(\xi_{\alpha,\beta}>z)dz\Bigr\}\\
&=q_{\alpha,\beta}(\lambda_0)\exp\Bigl\{-\int_0^\infty\frac{e^{-\lambda_0 z}-e^{-\lambda z}}{z}\mathbf{P}(\xi_{\alpha,\beta}>z)dz\Bigr\}.
\end{align*}
It follows from the definition of $\xi_{\alpha,\beta}$ that
\begin{equation}
\label{n7}
\mathbf{P}(\xi_{\alpha,\beta}>z)=\mathbf{P}\left(Y_{\alpha,\beta}(1)>z^{1-1/\alpha}\right)\sim\frac{C}{z^{\alpha-1}}
\text{ as }z\to\infty.
\end{equation}
This relation yields that
$$
\int_1^\infty\frac{1}{z}\mathbf{P}(\xi_{\alpha,\beta}>z)dz<\infty.
$$
Therefore,
\begin{align*}
\int_0^\infty\frac{e^{-\lambda_0 z}-e^{-\lambda z}}{z}\mathbf{P}(\xi_{\alpha,\beta}>z)dz\hspace{5cm}\\
=
\int_0^\infty\frac{1-e^{-\lambda z}}{z}\mathbf{P}(\xi_{\alpha,\beta}>z)dz
-\int_0^\infty\frac{1-e^{-\lambda_0 z}}{z}\mathbf{P}(\xi_{\alpha,\beta}>z)dz.
\end{align*}
Consequently,
$$
q_{\alpha,\beta}(\lambda)=C\exp\Bigl\{-\int_0^\infty\frac{1-e^{-\lambda z}}{z}\mathbf{P}(\xi_{\alpha,\beta}>z)dz\Bigr\}.
$$
To complete the proof of the theorem it remains to note that, in view of the scaling property of $Y_{\alpha,\beta}$,
$$
\mathbf{P}(\xi_{\alpha,\beta}>z)=\mathbf{P}\left(Y_{\alpha,\beta}(1)>z^{1-1/\alpha}\right)=
\mathbf{P}\left(Y_{\alpha,\beta}(z)-z>0\right).
$$
\section{Proof of Theorem \ref{T.exp}}\label{s_T.exp}
For every $\varepsilon\in(0,1)$,
\begin{equation}
\label{b1}
\mathbf{E}\bigl(\tau^{(a)}\bigr)^r=\sum_{n=0}^\infty[(n+1)^r-n^r]\mathbf{P}(\tau^{(a)}>n)=
\Sigma_1+\Sigma_2+\Sigma_3,
\end{equation}
where  
\begin{align*}
\Sigma_1&:=\sum_{0\leq n\leq\varepsilon n_a}[(n+1)^r-n^r]\mathbf{P}(\tau^{(a)}>n),\\
\Sigma_2&:=\sum_{\varepsilon n_a<n<n_a/\varepsilon}[(n+1)^r-n^r]\mathbf{P}(\tau^{(a)}>n),\\
\Sigma_3&:=\sum_{n\geq n_a/\varepsilon}[(n+1)^r-n^r]\mathbf{P}(\tau^{(a)}>n).
\end{align*}
Since $[(n+1)^r-n^r]\leq Cn^{r-1}$,
\begin{equation}
\label{b0}
\Sigma_1\leq
C\sum_{0\leq n\leq\varepsilon n_a}n^{r-1}\mathbf{P}(\tau^{(0)}>n)\leq C \varepsilon^{\rho+r-1} n_a^{r}\mathbf{P}(\tau^{(0)}>n_a).
\end{equation}
In the last step we used the fact that $\mathbf{P}(\tau^{(0)}>n)$ is regularly varying
with index $\rho-1$. 

Furthermore, in view of (\ref{Invar}),
\begin{align*}
\nonumber
\psi_a(r;x)&:=\bigl(([xn_a]+1)^r-([xn_a])^r\bigr)\frac{\mathbf{P}(\tau^{(a)}>[xn_a])}{n_a^{r-1}\mathbf{P}(\tau^{(0)}>n_a)}\\
&=\frac{([xn_a]+1)^r-([xn_a])^r}{n_a^{r-1}}\frac{\mathbf{P}(\tau^{(a)}>[xn_a])}{\mathbf{P}(\tau^{(0)}>[xn_a])}
\frac{\mathbf{P}(\tau^{(0)}>[xn_a])}{\mathbf{P}(\tau^{(0)}>n_a)}\\
&\to rx^{r-1}(1-F_{\alpha,\beta}(x^{1-1/\alpha}))x^{\rho-1}
\quad\text{as }a\to0.
\end{align*}
Then, by dominated convergence,
\begin{align}
\label{b2}
\nonumber
\lim_{a\to0}\frac{\Sigma_2}{n_a^r\mathbf{P}(\tau^{(0)}>n_a)}&=\lim_{a\to0}
\int_\varepsilon^{1/\varepsilon}\psi_a(r;x)dx\\
&=\int_\varepsilon^{1/\varepsilon}x^{r-1}(1-F_{\alpha,\beta}(x^{1-1/\alpha}))x^{\rho-1}dx.
\end{align}
In view of Proposition~\ref{UpBound},
\begin{equation*}
\Sigma_3\leq 
C\mathbf{E}\tau^{(a)}\sum_{n\geq n_a/\varepsilon} n^{r-1}\frac{V(na)}{(na)^2}.
\end{equation*}
Since $V(x)$ varies regularly,
\begin{align*}
\sum_{n\geq n_a/\varepsilon} n^{r-1}\frac{V(na)}{(na)^2}&\sim a^{-r}\int_{an_a/\varepsilon}^\infty x^{r-3}V(x)dx\\
&\sim(\alpha-r)^{-1}\varepsilon^{\alpha-r}a^{-r}(an_a)^{r-2}V(an_a)\\
&\sim (\alpha-r)^{-1}\varepsilon^{\alpha-r}n_a^{r}\frac{V(an_a)}{(an_a)}\sim(\alpha-r)^{-1}\varepsilon^{\alpha-r}n_a^{r-1}.
\end{align*}
Here we used the relations 
$$
an_a\sim c_{n_a}\quad\text{as }a\to0
$$ 
and 
$$c_n^{-2}V(c_n)\sim n^{-1}\quad\text{as }n\to\infty.
$$
Consequently,
\begin{equation}
\label{b3}
\Sigma_3\leq C\varepsilon^{\alpha-r}\mathbf{E}\tau^{(a)}n_a^{r-1}.
\end{equation}
Substituting (\ref{b0})--(\ref{b3}) with $r=1$ into (\ref{b1}) with $r=1$, we have
\begin{equation*}
\limsup_{a\to0}\frac{\mathbf{E}\tau^{(a)}}{n_a\mathbf{P}(\tau^{(0)}>n_a)}
\leq\frac{1}{1-C\varepsilon^{\alpha-1}}\Bigl(\int_\varepsilon^{1/\varepsilon}(1-F_{\alpha,\beta}(x^{1-1/\alpha}))x^{\rho-1}dx
+C\varepsilon^\rho\Bigr).
\end{equation*}
Thus,
\begin{equation*}
\mathbf{E}\tau^{(a)}\leq Cn_a\mathbf{P}(\tau^{(0)}>n_a).
\end{equation*}
Applying this inequality to (\ref{b3}), we get
\begin{equation}
\label{b4}
\Sigma_3\leq C\varepsilon^{\alpha-r}n_a^{r}\mathbf{P}(\tau^{(0)}>n_a).
\end{equation}
Combining (\ref{b1}), (\ref{b0}), (\ref{b2}) and (\ref{b4}), we obtain
\begin{align*}
\liminf_{a\to0}\frac{\mathbf{E}\bigl(\tau^{(a)}\bigr)^r}{n_a^r\mathbf{P}(\tau^{(0)}>n_a)}
\geq\int_\varepsilon^{1/\varepsilon}x^{r-1}(1-F_{\alpha,\beta}(x^{1-1/\alpha}))x^{\rho-1}dx
\end{align*}
and
\begin{align*}
\limsup_{a\to0}\frac{\mathbf{E}\bigl(\tau^{(a)}\bigr)^r}{n_a^r\mathbf{P}(\tau^{(0)}>n_a)}\hspace{6cm}\\
\leq\int_\varepsilon^{1/\varepsilon}x^{r-1}(1-F_{\alpha,\beta}(x^{1-1/\alpha}))x^{\rho-1}dx
+C\varepsilon^{\rho+r-1}+C\varepsilon^{\alpha-r}.
\end{align*} 
The latter inequality yields
\begin{equation}
\label{b5}
\limsup_{a\to0}\frac{\mathbf{E}\bigl(\tau^{(a)}\bigr)^r}{n_a^r\mathbf{P}(\tau^{(0)}>n_a)}
<\infty.
\end{equation}
Hence, letting $\varepsilon\to0$,
\begin{equation}
\label{b6}
\lim_{a\to0}\frac{\mathbf{E}\bigl(\tau^{(a)}\bigr)^r}{n_a^r\mathbf{P}(\tau^{(0)}>n_a)}
=\int_0^{\infty}x^{r-1}(1-F_{\alpha,\beta}(x^{1-1/\alpha}))x^{\rho-1}dx.
\end{equation}
The integral $\int_0^{\infty}x^{r-1}(1-F_{\alpha,\beta}(x^{1-1/\alpha}))x^{\rho-1}dx$ is finite in view of (\ref{b5}).
Noting now that $n_a^r\mathbf{P}(\tau^{(0)}>n_a)$ is regularly varying with index $-\alpha(\rho+r-1)/(\alpha-1)$,
we complete the proof of the theorem.
\section{Proofs of large deviation results}
\subsection{Proof of Theorem~\ref{LD_stable}}
Since $an/c_n\to\infty$ there exists $N(n)$ satisfying
\begin{equation*}
\frac{aN(n)}{c_n}\to\infty\text{ and }N(n)=o(n).
\end{equation*}
We now split the sum in (\ref{Rec}) into two parts:
\begin{equation*}
\Sigma_1:=\sum_{k=0}^{N(n)}\mathbf{P}(\tau^{(a)}>k)\mathbf{P}(S_{n-k}^{(0)}>(n-k)a),
\end{equation*}
\begin{equation*}
\Sigma_2:=\sum_{k=N(n)+1}^{n-1}\mathbf{P}(\tau^{(a)}>k)\mathbf{P}(S_{n-k}^{(0)}>(n-k)a).
\end{equation*}

Since 
\begin{equation*}
\lim_{j\to\infty}\sup_{x>q_jc_j}\Big|\frac{\mathbf{P}(S_j^{(0)}>x)}{j\mathbf{P}(X>x)}-1\Big|=0
\end{equation*}
for any sequence $q_j\uparrow\infty$, we get the relation
\begin{align}
\label{T2.1}
\nonumber
\Sigma_1&=(1+o(1))n\mathbf{P}(X>na)\sum_{k=0}^{N(n)}\mathbf{P}(\tau^{(a)}>k)\\
&=(1+o(1))n\mathbf{P}(X>na)\Bigl(\mathbf{E}\tau^{(a)}-\sum_{k=N(n)+1}^{n-1}\mathbf{P}(\tau^{(a)}>k)\Bigr).
\end{align}
Noting that $N(n)\gg n_a$ and taking into account (\ref{b3}), we see that
\begin{equation}
\label{T2.2}
\sum_{k=N(n)+1}^{n-1}\mathbf{P}(\tau^{(a)}>k)=o\Bigl(\mathbf{E}\tau^{(a)}\Bigr).
\end{equation}
Combining (\ref{T2.1}) and (\ref{T2.2}), we have
\begin{equation}
\label{T2.3}
\Sigma_1=(1+o(1))n\mathbf{E}\tau^{(a)}\mathbf{P}(X>na).
\end{equation}

We now turn our attention to $\Sigma_2$. It follows from Proposition~\ref{UpBound} that
\begin{align*}
\Sigma_2&\leq\mathbf{P}(\tau^{(a)}>N(n))\sum_{j=1}^n\mathbf{P}(S_j^{(0)}>aj)\\
&\leq C\mathbf{E}\tau^{(a)}\frac{V(aN(n))}{(aN(n))^2}\sum_{j=1}^n\mathbf{P}(S_j^{(0)}>aj).
\end{align*}
Furthermore, using (\ref{P3}), we obtain
\begin{equation}
\label{T2.4}
\Sigma_2\leq\mathbf{E}\tau^{(a)}\frac{V(aN(n))}{(aN(n))^2}n^2\mathbf{P}(|X|\geq na).
\end{equation}
From the definition of $c_n$ and the relation $aN(n)\gg c_n$ we conclude that
\begin{equation*}
\frac{V(aN(n))}{(aN(n))^2}=o(1/n).
\end{equation*} 
Moreover, $\mathbf{P}(|X|\geq na)\leq C\mathbf{P}(X\geq na)$
for every $X\in\mathcal{D}(\alpha,\beta)$ with $\alpha<2$ and $\beta>-1$. Then, (\ref{T2.4}) implies
\begin{equation}
\label{T2.5}
\Sigma_2=o\Bigl(n\mathbf{E}\tau^{(a)}\mathbf{P}(X>na)\Bigr).
\end{equation}
Substituting (\ref{T2.3}) and (\ref{T2.5}) into (\ref{Rec}), we complete the proof.
\subsection{Proof of Theorem~\ref{LD_var}}
Recall definition (\ref{xi-def}) of $\xi(a)$.
Set
\begin{equation}
\label{phi}
\phi_j:=e^{\xi(a)j}\mathbf{P}(\tau^{(a)}>j)\quad\text{and}\quad\theta_j:=e^{\xi(a)j}\mathbf{P}(S_j^{(0)}>aj).
\end{equation}
It is easily seen that
\begin{equation}
\label{3.0}
\theta_j\leq C\quad\text{for all }j\leq 1/a^2.
\end{equation}
Furthermore, combining (\ref{cond1}) with the relations
\begin{equation*}
\overline{\Phi}(x)\leq \frac{1}{x\sqrt{2\pi}}e^{-x^2/2}
\end{equation*}
and 
\begin{equation*}
\overline{\Phi}(x)\sim \frac{1}{x\sqrt{2\pi}}e^{-x^2/2}\quad\text{as }x\to\infty,
\end{equation*}
we get
\begin{equation}
\label{3.1}
\theta_j\leq\frac{C}{a\sqrt{j}},\ j\leq n,
\end{equation}
and
\begin{equation}
\label{3.2}
\theta_j\sim\frac{1}{a\sqrt{2\pi j}}\quad\text{for }j\leq n\text{ and }ja^2\to\infty,
\end{equation}
respectively.

Multiplying both sides of (\ref{Rec}) by $e^{a^2n/2}$, we see that the sequence $\phi_j$
satisfies the equation
\begin{equation}
\label{3.3}
k\phi_k=\sum_{j=0}^{k-1}\phi_j\theta_{k-j},\quad k\geq1.
\end{equation}
If $n$ satisfies the conditions of the theorem, then, using (\ref{3.0}) and (\ref{3.1}), we have
$$
\sup_{n\geq1}\max_{j\leq n}\theta_j<\infty.
$$ 
Consequently,
\begin{equation*}
\phi_k\leq\frac{C}{k}\sum_{j=o}^k\phi_j\leq\frac{C}{k}\sum_{j=o}^n\phi_j
\end{equation*}
for all $k\leq n$. Setting $\sigma_n:=\sum_{j=0}^n\phi_j,$
we rewrite the latter bound as 
$$
\phi_k\leq\frac{C}{k}\sigma_n,\quad k\leq n.
$$
Now, applying this bound and (\ref{3.1}) to the terms on the right hand side
of (\ref{3.3}), we obtain for all $k\leq n$ the bound
\begin{align}
\label{3.4}
\nonumber
\phi_k&=\frac{1}{k}\sum_{0\leq j< k/2}\phi_j\theta_{k-j}+\frac{1}{k}\sum_{k/2\leq j<k}\phi_j\theta_{k-j}\\
\nonumber
&\leq\frac{C}{ak^{3/2}}\sum_{0\leq j< k/2}\phi_j+\frac{C\sigma_n}{k^2}\sum_{k/2\leq j<k}\theta_{k-j}\\
&\leq\frac{C\sigma_n}{ak^{3/2}}+\frac{C\sigma_n}{k^2}\sum_{1\leq j<k}\frac{1}{a\sqrt{j}}\leq\frac{C\sigma_n}{ak^{3/2}}.
\end{align}
This inequality allows us to determine the asymptotic behaviour of $\phi_n$.
First of all we note that (\ref{3.2}) yields
\begin{equation*}
\sum_{0\leq j\leq N(n)}\phi_j\theta_{n-j}\sim\frac{1}{a\sqrt{2\pi n}}\sum_{0\leq j\leq N(n)}\phi_j\quad\text{as }a\to0,
\end{equation*}
for every $N(n)=o(n)$. Moreover, by (\ref{3.4}),
\begin{equation}
\label{3.4'}
0\leq\sigma_n-\sum_{0\leq j\leq N(n)}\phi_j=\sum_{N(n)<j\leq n}\phi_j\leq\frac{C\sigma_n}{aN(n)}.
\end{equation}
Therefore, choosing $N(n)$ satisfying
\begin{equation*}
N(n)=o(n)\quad\text{and}\quad aN^2(n)\to\infty,
\end{equation*}
we have, as $a\to0$,
\begin{equation}
\label{3.5}
\sum_{0\leq j\leq N(n)}\phi_j\theta_{n-j}\sim\frac{\sigma_n}{a\sqrt{2\pi n}}.
\end{equation}
Further, it follows from (\ref{3.1}) and (\ref{3.4}) that
\begin{align}
\label{3.6}
\nonumber
\sum_{N(n)<j<n/2}\phi_j\theta_{n-j}&\leq\frac{C}{a\sqrt{n}}\sum_{N(n)<j<n/2}\phi_j\\
&\leq \frac{C}{a\sqrt{n}}\sum_{N(n)<j<n/2}\frac{\sigma_n}{aj^{3/2}}\leq
\frac{C\sigma_n}{a^2\sqrt{nN(n)}}
\end{align}
and
\begin{align}
\label{3.7}
\sum_{n/2\leq j<n}\phi_j\theta_{n-j}&\leq\frac{C\sigma_n}{an^{3/2}}\sum_{j=1}^n\theta_j
\leq\frac{C\sigma_n}{an^{3/2}}\sum_{j=1}^n\frac{1}{a\sqrt{j}}\leq\frac{C\sigma_n}{a^2n}.
\end{align}
Combining (\ref{3.5}) -- (\ref{3.7}) and recalling that $a^2N(n)\to\infty$, we get
\begin{equation*}
\sum_{j=0}^{n-1}\phi_j\theta_{n-j}\sim\frac{\sigma_n}{a\sqrt{2\pi n}}\quad\text{as }a\to0.
\end{equation*} 
Substituting this into (\ref{3.3}), we have
\begin{equation}
\label{3.8}
\phi_n\sim\frac{\sigma_n}{a\sqrt{2\pi}n^{3/2}}\quad\text{as }a\to0.
\end{equation}

To complete the proof of the theorem it remains to find the asymptotic behaviour of $\sigma_n$.
First of all, (\ref{3.4'}) implies that the bounds
\begin{equation}
\label{3.9}
\sum_{j\leq 1/\varepsilon a^2}\phi_j\leq\sigma_n\leq
(1-C\sqrt{\varepsilon})^{-1}\sum_{j\leq 1/\varepsilon a^2}\phi_j
\end{equation}
are valid for all sufficiently small values of $\varepsilon$.
Applying Theorem~\ref{T.normal} and recalling that $\mathbf{P}(\tau^{(0)}>j)$ is regularly varying
with index $-1/2$, we see that
\begin{align*}
\lim_{a\to0}\frac{\phi_{[xa^{-2}]}}{\mathbf{P}(\tau^{(0)}>a^{-2})}
&=\lim_{a\to0}\frac{e^{[xa^{-2}]\xi(a)}\mathbf{P}(\tau^{(a)}>[xa^{-2}])}{\mathbf{P}(\tau^{(0)}>[xa^{-2}])}
\frac{\mathbf{P}(\tau^{(0)}>[xa^{-2}])}{\mathbf{P}(\tau^{(0)}>a^{-2})}\\
&=e^{x/2}(1-F_{2,0}(\sqrt{x}))\frac{1}{\sqrt{x}}
\end{align*}
for every $x>0$.
Thus, by dominated convergence,
\begin{align}
\label{3.10}
\nonumber
\lim_{a\to0}\frac{\sum\limits_{j\leq 1/\varepsilon a^2}\phi_j}{a^{-2}\mathbf{P}(\tau^{(0)}>a^{-2})}
&=\int_0^{1/\varepsilon}\frac{e^{x/2}}{\sqrt{x}}(1-F_{2,0}(\sqrt{x}))dx\\
&=\int_0^{1/\varepsilon}\frac{e^{x/2}}{\sqrt{x}}(1-F_{2,0}(\sqrt{x}))dx
=:I(\varepsilon).
\end{align}
Using (\ref{F6}), we have
\begin{align*}
I(\varepsilon)&=\int_0^{1/\varepsilon}e^{x/2}\int_{\sqrt{x}}^\infty v^{-2}e^{-v^2/2}dvdx\\
&=\int_0^{\infty}e^{x/2}\int_{\sqrt{x}}^\infty v^{-2}e^{-v^2/2}dvdx-
\int_{1/\varepsilon}^\infty e^{x/2}\int_{\sqrt{x}}^\infty v^{-2}e^{-v^2/2}dvdx.
\end{align*}
Noting that
\begin{equation*}
\int_{\sqrt{x}}^\infty v^{-2}e^{-v^2/2}dv\leq\frac{e^{-x/2}}{x^{3/2}},
\end{equation*}
we have
\begin{equation*}
0\leq\int_0^{\infty}e^{x/2}\int_{\sqrt{x}}^\infty v^{-2}e^{-v^2/2}dvdx-I(\varepsilon)
\leq\sqrt{\varepsilon}.
\end{equation*}
Changing the order of integration and substituting $v^2/2=u$, we have
\begin{align*}
\int_0^{\infty}e^{x/2}\int_{\sqrt{x}}^\infty v^{-2}e^{-v^2/2}dvdx=
\int_0^\infty v^{-2}e^{-v^2/2}\int_0^{v^2}e^{x/2}dxdv\hspace{1cm}\\
=2\int_0^\infty v^{-2}e^{-v^2/2}(1-e^{-v^2/2})dv=\frac{1}{\sqrt{2}}\int_0^\infty u^{-3/2}(1-e^{-u})du.
\end{align*}
Integrating now by parts, we get
\begin{equation*}
\int_0^\infty u^{-3/2}(1-e^{-u})du=2\int_0^\infty u^{-1/2}e^{-u}du=2\Gamma(1/2)=2\sqrt{\pi}.
\end{equation*}
As a result we have the bounds
\begin{equation}
\label{3.11}
\sqrt{2\pi}-\sqrt{\varepsilon}\leq I(\varepsilon)\leq\sqrt{2\pi}.
\end{equation}
Substituting (\ref{3.10}) and (\ref{3.11}) into (\ref{3.9}), we obtain
\begin{align*}
\sqrt{2\pi}-\sqrt{\varepsilon}&\leq\liminf_{a\to0}\frac{\sigma_n}{a^{-2}\mathbf{P}(\tau^{(0)}>a^{-2})}\\
&\leq\limsup_{a\to0}\frac{\sigma_n}{a^{-2}\mathbf{P}(\tau^{(0)}>a^{-2})}\leq
\frac{\sqrt{2\pi}}{1-C\sqrt{\varepsilon}}.
\end{align*}
Since $\varepsilon$ can be chosen arbitrary small,
\begin{equation}
\label{3.12}
\sigma_n\sim\sqrt{2\pi}a^{-2}\mathbf{P}(\tau^{(0)}>a^{-2})
\end{equation} 
Combining (\ref{3.8}) and (\ref{3.12}), and recalling definition (\ref{phi}) of $\phi_n$, we have
\begin{equation}
\label{3.13}
\mathbf{P}(\tau^{(a)}>n)\sim a^{-3}n^{-3/2}e^{-\xi(a)n}\mathbf{P}(\tau^{(0)}>a^{-2}).
\end{equation}
Further, it follows from (\ref{b6}) that
$$
\mathbf{E}\tau^{(a)}\sim a^{-2}\mathbf{P}(\tau^{(0)}>a^{-2})\int_0^\infty (1-F_{2,0}(\sqrt{x}))x^{-1/2}dx.
$$
Substituting $\sqrt{x}=y$ and using (\ref{F6}), we get
\begin{align*}
\int_0^\infty (1-F_{2,0}(\sqrt{x}))x^{-1/2}dx=2\int_0^\infty (1-F_{2,0}(y))dy\\
=2\int_0^\infty y\Bigl(\int_y^\infty v^{-2}e^{-v^2/2}dv\Bigr)dy
=2\int_0^\infty v^{-2}e^{-v^2/2}\Bigl(\int_0^v ydy\Bigr)dv\\
=\int_0^\infty e^{-v^2/2}dv=\sqrt{\frac{\pi}{2}}.
\end{align*}
As a result we have
\begin{equation}
\label{3.14}
a^{-2}\mathbf{P}(\tau^{(0)}>a^{-2})\sim \sqrt{\frac{2}{\pi}}\mathbf{E}\tau^{(a)}.
\end{equation}
Combining (\ref{3.13}) and (\ref{3.14}), and noting that
\begin{equation*}
\frac{1}{a\sqrt{2\pi n}}e^{-\xi(a)n}\sim\overline{\Phi}(a\sqrt{n})
\exp\{na^3\lambda_m(a)\},
\end{equation*}
we complete the proof of the theorem.
\subsection{Proof of Theorem~\ref{tail}}
It is easy to see that there exist a constant $C$ 
and a regularly varying function $N(a)$ such that
\begin{equation}
\label{5.0}
\lim_{a\to0}\frac{N(a)}{a^{-2}\log a^{-2}}=(p-2)^{1/2}
\end{equation}
and
\begin{equation}
\label{5.1}
\sup_{n\leq N(a)}\frac{nP(X\geq na+\sqrt{n})}{\overline{\Phi}(a\sqrt{n})}\leq C
\text{ and }
\sup_{n\geq N(a)}\frac{\overline{\Phi}(a\sqrt{n})}{nP(X\geq na+\sqrt{n})}\leq C.
\end{equation}
We now split the right hand side of (\ref{GenF}) into the product of two exponentials:
\begin{align*}
&\exp\left\{\sum_{n=1}^\infty\frac{z^n}{n}\mathbf{P}(S_n^{(a)}>0)\right\}\\
&=\exp\left\{\sum_{n=1}^{N(a)}\frac{z^n}{n}\mathbf{P}(S_n^{(a)}>0)\right\}
\exp\left\{\sum_{n=N(a)+1}^\infty\frac{z^n}{n}\mathbf{P}(S_n^{(a)}>0)\right\}\\
&=:\left(\sum_{n=0}^\infty \psi_{1,n}z^n\right)\left(1+\sum_{n=N(a)+1}^\infty \psi_{2,n}z^n\right).
\end{align*}
Therefore,
\begin{equation}
\label{5.2}
\mathbf{P}(\tau^{(a)}>n)=\psi_{1,n}+\sum_{k=N(a)+1}^n\psi_{1,n-k}\psi_{2,k},\ n\geq1.
\end{equation}

We first want to find the asymptotic behaviour of $\psi_{2,n}$. We start by noting that
\begin{equation}
\label{5.3}
\psi_{2,n}=\sum_{j=1}^\infty\frac{1}{j!}q^{*j}_n,\ n> N(a),
\end{equation}
where $\{q^{*j}_n,n\geq1\}$ is the $j$-th convolution of 
$\displaystyle\left\{n^{-1}\mathbf{P}(S_n^{(a)}>0){\rm 1}\{n>N(a)\},n\geq1\right\}$. 
It follows from the second inequality in (\ref{5.1}) that
\begin{align*}
q^{*2}_n&=\sum_{k=N(a)+1}^{n-N-1}\frac{1}{k}\mathbf{P}(S_k^{(a)}>0)\frac{1}{n-k}\mathbf{P}(S_{n-k}^{(a)}>0)\\
&\leq C\sum_{k=N(a)+1}^{n-N-1}\mathbf{P}(X\geq ak)\mathbf{P}(X\geq a(n-k))\\
&\leq C\mathbf{P}(X\geq an/2)\sum_{N(a)+1}\mathbf{P}(X\geq ak)
\leq C\mathbf{P}(X\geq an)\int_{N(a)}^\infty\mathbf{P}(X\geq ay)dy.
\end{align*}

Since $\mathbf{P}(X\geq y)$ is regularly varying, we have
$$
\int_{N(a)}^\infty\mathbf{P}(X\geq ay)dy=
\frac{1}{a}\int_{aN(a)}^\infty\mathbf{P}(X\geq y)dy\leq
CN(a)\mathbf{P}(X\geq aN(a)).
$$
From this bound and (\ref{5.0}) we get
$$
q^{*2}_n\leq G(a)\mathbf{P}(X\geq an),
$$
where $G$ is regularly varying with index $p-2>0$. Then, by induction,
\begin{equation}
\label{5.4}
q^{*j}_n\leq G(a)\mathbf{P}(X\geq an)\text{ for all }j\geq2.
\end{equation}
Combining (\ref{5.3}), (\ref{5.4}), and using (\ref{NagRoz}) and (\ref{5.1}), we obtain the bound
\begin{align}
\nonumber
\psi_{2,n}&=\mathbf{P}(S_n^{(a)}>0)+\sum_{j=2}^\infty q^{*j}_n\\
&\leq C\Bigl(\frac{1}{n}\overline{\Phi}(a\sqrt{n})+\mathbf{P}(X\geq an)+G(a)\mathbf{P}(X\geq an)\Bigr)
\leq C\mathbf{P}(X\geq an)
\label{5.5}
\end{align}
and, for $n\geq Ca^{-2}\log a^{-2}$ with some $C>(p-2)^{1/2}$, the relation
\begin{align}
\label{5.6}
\nonumber
\psi_{2,n}&=\mathbf{P}(S_n^{(a)}>0)+O\Bigl(G(a)\mathbf{P}(X\geq an)\Bigr)\\
&\sim\frac{1}{n}\overline{\Phi}(a\sqrt{n})+\mathbf{P}(X\geq an)\sim \mathbf{P}(X\geq an).
\end{align}
In the last step we have used the fact that 
$\overline{\Phi}(a\sqrt{n})=o(\mathbf{P}(X\geq an))$ for $n\geq Ca^{-2}\log a^{-2}$, $C>(p-2)^{1/2}$.

From the first inequality in (\ref{5.1}) and (\ref{NagRoz}), which is valid under the condition (\ref{Roz_cond}),
we conclude that
$$
\mathbf{P}(S_n^{(a)}>0)\leq C\overline{\Phi}(a\sqrt{n})
$$
for all $n\leq N(a)$. Using arguments from the proof of Theorem~\ref{LD_var}, one sees that
\begin{equation}
\label{5.7}
\psi_{1,k}\leq\frac{C}{k}\overline{\Phi}(a\sqrt{k}),\ k\geq1.
\end{equation}
Combining (\ref{5.5}) and (\ref{5.7}), and applying the second inequality in (\ref{5.1}), we get
\begin{align*}
\sum_{k=N(a)}^{n-N(a)}\psi_{1,n-k}\psi_{2,k}
&\leq C\sum_{k=N(a)}^{n-N(a)}\frac{1}{n-k}\overline{\Phi}(a\sqrt{n-k})\mathbf{P}(X\geq ak)\\
&\leq C\sum_{k=N(a)}^{n-N(a)}\mathbf{P}\bigl(X\geq a(n-k)\bigr)\mathbf{P}(X\geq ak).
\end{align*}
In the derivation of (\ref{5.4}) we have shown that the sum in the last line is bounded by
$G(a)\mathbf{P}(X\geq an)$. Hence,
\begin{equation}
\label{5.8}\sum_{k=N(a)}^{n-N(a)}\psi_{1,n-k}\psi_{2,k}=O\Bigl(G(a)\mathbf{P}(X\geq an)\Bigr).
\end{equation}

It follows  from (\ref{GenF}) and the definition of $\{\psi_{1,n},n\geq1\}$ that
$\psi_{1,k}=\mathbf{P}(\tau^{(a)}>k)$ for all $k\leq N(a)$. Consequently,
\begin{align*}
\sum_{k=n-N(a)+1}^n\psi_{1,n-k}\psi_{2,k}&=
\sum_{k=0}^{N(a)-1}\mathbf{P}(\tau^{(a)}>k)\psi_{2,n-k}\\
&=\sum_{k=0}^{\tilde{N}(a)-1}\mathbf{P}(\tau^{(a)}>k)\psi_{2,n-k}+
\sum_{k=\tilde{N}(a)}^{N(a)-1}\mathbf{P}(\tau^{(a)}>k)\psi_{2,n-k},
\end{align*}
where $\tilde{N}(a)$ is such that $a^{-2}\ll \tilde{N}(a)\ll a^{-2}\log a^{-2}$.
Applying (\ref{5.5}) to the fist sum and (\ref{5.6}) to the second sum, we get 
\begin{align*}
\sum_{k=n-N(a)+1}^n\psi_{1,n-k}\psi_{2,k}=(1+o(1))\mathbf{P}(X\geq an)\sum_{k=0}^{\tilde{N}(a)-1}\mathbf{P}(\tau^{(a)}>k)\\
+O\Bigl(\mathbf{P}(X\geq an)\sum_{k=\tilde{N}(a)}^{N(a)-1}\mathbf{P}(\tau^{(a)}>k)\Bigr).
\end{align*}
Note that (\ref{b3}) implies that 
$$
\sum_{k=\tilde{N}(a)}^\infty\mathbf{P}(\tau^{(a)}>k)=o(\mathbf{E}\tau^{(a)}).
$$
Hence we finally obtain
\begin{equation}
\label{5.9}
\sum_{k=n-N(a)+1}^n\psi_{1,n-k}\psi_{2,k}\sim \mathbf{E}\tau^{(a)}\mathbf{P}(X\geq an).
\end{equation}
Combining (\ref{5.2}), (\ref{5.8}) and (\ref{5.9}), we have
$$
\mathbf{P}(\tau^{(a)}>n)=(1+o(1))\mathbf{E}\tau^{(a)}\mathbf{P}(X\geq an)+\psi_{1,n}.
$$
In order to finish the proof it remains to apply (\ref{5.7}) and to note that
$n^{-1}\overline{\Phi}(a\sqrt{n})=o(\mathbf{P}(X\geq an))$.

\vspace{12pt}

{\em Acknowledgements.} The author would like to thank Klaus Fleischmann for remarks which led to
a better presentation of the results.

\end{document}